\documentclass{amsart}
\usepackage{amsthm}
\usepackage{amstext}
\usepackage{amssymb}
\usepackage{hyperref}

\makeatletter
\numberwithin{equation}{section}
\numberwithin{figure}{section}
\theoremstyle{plain}
\newtheorem{thm}{Theorem}[section]
  \theoremstyle{remark}
  \newtheorem*{rem*}{Remark}

\allowdisplaybreaks
\usepackage[sort&compress,square,numbers]{natbib}

\newcommand\relphantom[1]{\mathrel{\phantom{#1}}} 
\global\long\def\Li{\mathrm{Li}}

\makeatother

\begin{document}

\title{A note on degenerate poly-Bernoulli numbers and polynomials}

\author{Dae San Kim}
\address{Department of Mathematics, Sogang University, Seoul 121-742, Republic
of Korea}
\email{dskim@sogang.ac.kr}

\author{Taekyun Kim}
\address{Department of Mathematics, Kwangwoon University, Seoul 139-701, Republic
of Korea}
\email{tkkim@kw.ac.kr}


\begin{abstract}
In this paper, we consider the degenerate poly-Bernoulli polynomials
and present new and explicit formulas for computing them in terms of the degenerate Bernoulli polynomials and Stirling numbers of the second kind.
\end{abstract}

\maketitle

\section{Introduction}

For $\lambda\in\mathbb{C}$, L. Carlitz considered the degenerate
Bernoulli polynomials given by the generating function 
\begin{equation}
\frac{t}{\left(1+\lambda t\right)^{\frac{1}{\lambda}}-1}\left(1+\lambda t\right)^{\frac{x}{\lambda}}=\sum_{n=0}^{\infty}\beta_{n}\left(x\mid\lambda\right)\frac{t^{n}}{n!},\quad\left(\text{see \cite{key-6,key-12,key-19}}\right).\label{eq:1}
\end{equation}

When $x=0$, $\beta_{n}\left(\lambda\right)=\beta_{n}\left(0\mid\lambda\right)$
are called the degenerate Bernoulli numbers. 

Thus, by (\ref{eq:1}), we get 
\begin{equation}
\beta_{n}\left(x\mid\lambda\right)=\sum_{l=0}^{n}\binom{n}{l}\beta_{l}\left(\lambda\right)\left(x\mid\lambda\right)_{n-l},\label{eq:2}
\end{equation}
where $\left(x\mid\lambda\right)_{n}=x\left(x-\lambda\right)\left(x-2\lambda\right)\cdots\left(x-\lambda\left(n-1\right)\right)$. 

The classical polylograithm function $\Li_{k}$ is 
\begin{equation}
\Li_{k}\left(x\right)=\sum_{n=1}^{\infty}\frac{x^{n}}{n^{k}},\quad\left(k\in\mathbb{Z}\right),\quad\left(\text{see \cite{key-10,key-11,key-12,key-13,key-14,key-15,key-16,key-17,key-18}}\right).\label{eq:3}
\end{equation}

From (\ref{eq:1}), we note that 
\begin{align}
 & \relphantom{=}{}\sum_{n=0}^{\infty}\lim_{\lambda\rightarrow0}\beta_{n}\left(x\mid\lambda\right)\frac{t^{n}}{n!}\label{eq:4}\\
 & =\lim_{\lambda\rightarrow0}\frac{t}{\left(1+\lambda t\right)^{\frac{1}{\lambda}}-1}\left(1+\lambda t\right)^{\frac{x}{\lambda}}\nonumber \\
 & =\frac{t}{e^{t}-1}e^{xt}\nonumber \\
 & =\sum_{n=0}^{\infty}B_{n}\left(x\right)\frac{t^{n}}{n!},\nonumber 
\end{align}
where $B_{n}\left(x\right)$ are called the Bernoulli polynomials
(see \cite{key-1,key-2,key-3,key-4,key-5,key-6,key-7,key-8,key-9,key-10,key-11,key-12,key-13,key-14,key-15,key-16,key-17,key-18,key-19,key-20,key-21,key-22,key-23,key-24,key-25,key-26,key-27}). 

Thus, by (\ref{eq:4}), we get 
\begin{equation}
\lim_{\lambda\rightarrow0}\beta_{n}\left(x\mid\lambda\right)=B_{n}\left(x\right),\quad\left(n\ge0\right).\label{eq:5}
\end{equation}

In \cite{key-3,key-10}, the poly-Bernoulli polynomials are given
by 
\begin{equation}
\frac{\Li_{k}\left(1-e^{-t}\right)}{e^{t}-1}e^{xt}=\sum_{n=0}^{\infty}B_{n}^{\left(k\right)}\left(x\right)\frac{t^{n}}{n!}.\label{eq:6}
\end{equation}

For $k=1$, we have 
\begin{equation}
\frac{\Li_{1}\left(1-e^{-t}\right)}{e^{t}-1}e^{xt}=\frac{t}{e^{t}-1}e^{xt}=\sum_{n=0}^{\infty}B_{n}\left(x\right)\frac{t^{n}}{n!}.\label{eq:7}
\end{equation}

By (\ref{eq:4}) and (\ref{eq:7}), we get $B_{n}^{\left(1\right)}\left(x\right)=B_{n}\left(x\right)$. 

The Stirling numbers of the second kind are given by 
\begin{equation}
x^{n}=\sum_{l=0}^{n}S_{2}\left(n,l\right)\left(x\right)_{l},\quad\left(\text{see \cite{key-1,key-2,key-3,key-4,key-5,key-6,key-7,key-8,key-9,key-10,key-11,key-12,key-13,key-14,key-15,key-16,key-17,key-18,key-19,key-20,key-21,key-22,key-23,key-24,key-25,key-26,key-27}}\right).\label{eq:8}
\end{equation}
and the Stirling numbers of the first kind are defined by 
\begin{equation}
\left(x\right)_{n}=x\left(x-1\right)\cdots\left(x-n+1\right)=\sum_{l=0}^{n}S_{1}\left(n,l\right)x^{l},\quad\left(n\ge0\right).\label{eq:9}
\end{equation}

The purpose of this paper is to construct the degenerate poly-Bernoulli
polynomials and present new and explicit formulas for computing them
in terms of the degenerate Bernoulli polynomials and Stirling
numbers of the second kind.

\section{Degenerate poly-Bernoulli numbers and polynomials}

For $\lambda\in\mathbb{C}$, $k\in\mathbb{Z}$, we consider the degenerate
poly-Bernoulli polynomials given by the generating function

\begin{equation}
\frac{\Li_{k}\left(1-e^{-t}\right)}{\left(1+\lambda t\right)^{\frac{1}{\lambda}}-1}\left(1+\lambda t\right)^{\frac{x}{\lambda}}=\sum_{n=0}^{\infty}\beta_{n}^{\left(k\right)}\left(x\mid\lambda\right)\frac{t^{n}}{n!}.\label{eq:10}
\end{equation}

When $x=0$, $\beta_{n}^{\left(k\right)}\left(\lambda\right)=\beta_{n}^{\left(k\right)}\left(0\mid\lambda\right)$
are called the degenerate poly-Bernoulli numbers. Note that $\beta_{n}^{\left(1\right)}\left(x\mid\lambda\right)=\beta_{n}\left(x\mid\lambda\right)$
and $\lim_{\lambda\rightarrow0}\beta_{n}^{\left(k\right)}\left(x\mid\lambda\right)=B_{n}^{\left(k\right)}\left(x\right)$. 

From (\ref{eq:10}), we can derive the following equation: 
\begin{align}
\sum_{n=0}^{\infty}\beta_{n}^{\left(k\right)}\left(x\mid\lambda\right)\frac{t^{n}}{n!} & =\left(\frac{\Li_{k}\left(1-e^{-t}\right)}{\left(1+\lambda t\right)^{\frac{1}{\lambda}}-1}\right)\left(1+\lambda t\right)^{\frac{x}{\lambda}}\label{eq:11}\\
 & =\left(\sum_{l=0}^{\infty}\beta_{l}^{\left(k\right)}\left(\lambda\right)\frac{t^{l}}{l!}\right)\left(\sum_{m=0}^{\infty}\left(x\mid\lambda\right)_{m}\frac{t^{m}}{m!}\right)\nonumber \\
 & =\sum_{n=0}^{\infty}\left(\sum_{l=0}^{n}\binom{n}{l}\beta_{l}^{\left(k\right)}\left(\lambda\right)\left(x\mid\lambda\right)_{n-l}\right)\frac{t^{n}}{n!}.\nonumber 
\end{align}

Thus, by (\ref{eq:11}), we get 
\begin{equation}
\beta_{n}^{\left(k\right)}\left(x\mid\lambda\right)=\sum_{l=0}^{n}\binom{n}{l}\beta_{l}^{\left(k\right)}\left(\lambda\right)\left(x\mid\lambda\right)_{n-l}.\label{eq:12}
\end{equation}

Now, we observe that 
\begin{align}
 & \relphantom{=}{}\frac{\Li_{k}\left(1-e^{-t}\right)}{\left(1+\lambda t\right)^{\frac{1}{\lambda}}-1}\left(1+t\right)^{\frac{x}{\lambda}}\label{eq:13}\\
 & =\sum_{n=0}^{\infty}\beta_{n}^{\left(k\right)}\left(x\mid\lambda\right)\frac{t^{n}}{n!}\nonumber \\
 & =\frac{\left(1+t\right)^{\frac{x}{\lambda}}}{\left(1+\lambda t\right)^{\frac{1}{\lambda}}-1}\underset{\left(k-2\right)\text{times}}{\int_{0}^{t}\underbrace{\frac{1}{e^{y}-1}\int_{0}^{y}\frac{1}{e^{y}-1}\int_{0}^{y}\cdots\frac{1}{e^{y}-1}\int_{0}^{y}}}\frac{y}{e^{y}-1}dy\cdots dy.\nonumber 
\end{align}

From (\ref{eq:13}), we have 
\begin{align}
 & \relphantom{=}{}\sum_{n=0}^{\infty}\beta_{n}^{\left(2\right)}\left(x\mid\lambda\right)\frac{t^{n}}{n!}\label{eq:14}\\
 & =\frac{\left(1+t\right)^{\frac{x}{\lambda}}}{\left(1+\lambda t\right)^{\frac{1}{\lambda}}-1}\int_{0}^{t}\frac{y}{e^{y}-1}dy\nonumber \\
 & =\frac{\left(1+t\right)^{\frac{x}{\lambda}}}{\left(1+\lambda t\right)^{\frac{1}{\lambda}}-1}\sum_{l=0}^{\infty}\frac{B_{l}}{l!}\int_{0}^{t}y^{l}dy\nonumber \\
 & =\left(\frac{t}{\left(1+\lambda t\right)^{\frac{1}{\lambda}}-1}\left(1+\lambda t\right)^{\frac{x}{\lambda}}\right)\left(\sum_{l=0}^{\infty}\frac{B_{l}}{l+1}\frac{t^{l}}{l!}\right)\nonumber \\
 & =\sum_{n=0}^{\infty}\left\{ \sum_{l=0}^{n}\binom{n}{l}\frac{B_{l}}{l+1}\beta_{n-l}\left(x\mid\lambda\right)\right\} \frac{t^{n}}{n!},\nonumber 
\end{align}
where $B_{n}=B_{n}\left(0\right)$ are Bernoulli numbers.

By comparing the coefficients on both sides of (\ref{eq:14}),
we obtain the following theorem.
\begin{thm}
\label{thm:1} For $n\ge0$, we have 
\begin{align*}
\beta_{n}^{\left(2\right)}\left(x\mid\lambda\right) & =\sum_{l=0}^{n}\binom{n}{l}\frac{B_{l}}{l+1}\beta_{n-l}\left(x\mid\lambda\right)\\
 & =\beta_{n}\left(x\mid\lambda\right)-\frac{n}{4}\beta_{n-1}\left(x\mid\lambda\right)+\sum_{l=2}^{n}\binom{n}{l}\frac{B_{l}}{l+1}\beta_{n-l}\left(x\mid\lambda\right).
\end{align*}

Moreover, 
\[
\beta_{n}^{\left(k\right)}\left(x\mid\lambda\right)=\sum_{l=0}^{n}\binom{n}{l}\beta_{l}^{\left(k\right)}\left(\lambda\right)\left(x\mid\lambda\right)_{n-l}.
\]

\end{thm}
By (\ref{eq:13}), we easily get
\begin{align}
 & \relphantom{=}{}\sum_{n=0}^{\infty}\beta_{n}^{\left(k\right)}\left(x\mid\lambda\right)\frac{t^{n}}{n!}\label{eq:15}\\
 & =\frac{\Li_{k}\left(1-e^{-t}\right)}{\left(1+\lambda t\right)^{\frac{1}{\lambda}}-1}\left(1+t\right)^{\frac{x}{\lambda}}\nonumber \\
 & =\frac{t}{\left(1+\lambda t\right)^{\frac{1}{\lambda}}-1}\left(1+t\right)^{\frac{x}{\lambda}}\frac{\Li_{k}\left(1-e^{-t}\right)}{t}.\nonumber 
\end{align}

We observe that 

\begin{align}
\frac{1}{t}\Li_{k}\left(1-e^{-t}\right) & =\frac{1}{t}\sum_{n=1}^{\infty}\frac{1}{n^{k}}\left(1-e^{-t}\right)^{n}\label{eq:16}\\
 & =\frac{1}{t}\sum_{n=1}^{\infty}\frac{\left(-1\right)^{n}}{n^{k}}n!\sum_{l=n}^{\infty}S_{2}\left(l,n\right)\frac{\left(-t\right)^{l}}{l!}\nonumber \\
 & =\frac{1}{t}\sum_{l=1}^{\infty}\sum_{n=1}^{l}\frac{\left(-1\right)^{n+l}}{n^{k}}n!S_{2}\left(l,n\right)\frac{t^{l}}{l!}\nonumber \\
 & =\sum_{l=0}^{\infty}\sum_{n=1}^{l+1}\frac{\left(-1\right)^{n+l+1}}{n^{k}}n!\frac{S_{2}\left(l+1,n\right)}{l+1}\frac{t^{l}}{l!}.\nonumber 
\end{align}

From (\ref{eq:15}) and (\ref{eq:16}), we have 
\begin{align}
 & \relphantom{=}{}\sum_{n=0}^{\infty}\beta_{n}^{\left(k\right)}\left(x\mid\lambda\right)\frac{t^{n}}{n!}\label{eq:17}\\
 & =\left(\sum_{m=0}^{\infty}\beta_{m}\left(x\mid\lambda\right)\frac{t^{m}}{m!}\right)\left(\sum_{l=0}^{\infty}\left(\sum_{p=1}^{l+1}\frac{\left(-1\right)^{p+l+1}}{p^{k}}p!\frac{S_{2}\left(l+1,p\right)}{l+1}\right)\frac{t^{l}}{l!}\right)\nonumber \\
 & =\sum_{n=0}^{\infty}\left\{ \sum_{l=0}^{n}\binom{n}{l}\left(\sum_{p=1}^{l+1}\frac{\left(-1\right)^{p+l+1}p!}{p^{k}}\frac{S_{2}\left(l+1,p\right)}{l+1}\right)\beta_{n-l}\left(x\mid\lambda\right)\right\} \frac{t^{n}}{n!}.\nonumber 
\end{align}

By comparing the coefficients on both sides of (\ref{eq:17}),
we obtain the following theorem.
\begin{thm}
\label{thm:2} For $n\ge0$, we have 
\[
\beta_{n}^{\left(k\right)}\left(x\mid\lambda\right)=\sum_{l=0}^{n}\binom{n}{l}\left(\sum_{p=1}^{l+1}\frac{\left(-1\right)^{p+l+1}p!}{p^{k}}\frac{S_{2}\left(l+1,p\right)}{l+1}\right)\beta_{n-l}\left(x\mid\lambda\right).
\]

\end{thm}
It is easy to show that 
\begin{align}
 & \relphantom{=}\frac{\Li_{k}\left(1-e^{-t}\right)}{\left(1+\lambda t\right)^{\frac{1}{\lambda}}-1}\left(1+\lambda t\right)^{\frac{x+1}{\lambda}}-\frac{\Li_{k}\left(1-e^{-t}\right)}{\left(1+\lambda t\right)^{\frac{1}{\lambda}}-1}\left(1+\lambda t\right)^{\frac{x}{\lambda}}\label{eq:18}\\
 & =\left(1+\lambda t\right)^{\frac{x}{\lambda}}\Li_{k}\left(1-e^{-t}\right)\nonumber \\
 & =\left(\sum_{l=0}^{\infty}\left(x\mid\lambda\right)_{l}\frac{t^{l}}{l!}\right)\left(\sum_{m=1}^{\infty}\frac{\left(1-e^{-t}\right)^{m}}{m^{k}}\right)\nonumber \\
 & =\left(\sum_{l=0}^{\infty}\left(x\mid\lambda\right)_{l}\frac{t^{l}}{l!}\right)\left(\sum_{m=0}^{\infty}\frac{\left(1-e^{-t}\right)^{m+1}}{\left(m+1\right)^{k}}\right)\nonumber \\
 & =\left(\sum_{l=0}^{\infty}\left(x\mid\lambda\right)_{l}\frac{t^{l}}{l!}\right)\left(\sum_{p=1}^{\infty}\left(\sum_{m=0}^{p-1}\frac{\left(-1\right)^{m+p+1}}{\left(m+1\right)^{k}}\left(m+1\right)!S_{2}\left(p,m+1\right)\right)\frac{t^{p}}{p!}\right)\nonumber \\
 & =\sum_{n=1}^{\infty}\left\{ \sum_{p=1}^{n}\sum_{m=0}^{p-1}\frac{\left(-1\right)^{m+p+1}}{\left(m+1\right)^{k}}\left(m+1\right)!S_{2}\left(p,m+1\right)\binom{n}{p}\left(x\mid\lambda\right)_{n-p}\right\} \frac{t^{n}}{n!}.\nonumber 
\end{align}

On the other hand, 
\begin{align}
 & \relphantom{=}{}\frac{\Li_{k}\left(1-e^{-t}\right)}{\left(1+\lambda t\right)^{\frac{1}{\lambda}}-1}\left(1+\lambda t\right)^{\frac{x+1}{\lambda}}-\frac{\Li_{k}\left(1-e^{-t}\right)}{\left(1+\lambda t\right)^{\frac{1}{\lambda}}-1}\left(1+\lambda t\right)^{\frac{x}{\lambda}}\label{eq:19}\\
 & =\sum_{n=0}^{\infty}\left\{ \beta_{n}^{\left(k\right)}\left(x+1\mid\lambda\right)-\beta_{n}^{\left(k\right)}\left(x\mid\lambda\right)\right\} \frac{t^{n}}{n!}.\nonumber 
\end{align}

Therefore, by (\ref{eq:18}) and (\ref{eq:19}), we obtain the following
theorem.
\begin{thm}
\label{thm:3} For $n\ge1$, we have 
\begin{align*}
 & \relphantom{=}{}\beta_{n}^{\left(k\right)}\left(x+1\mid\lambda\right)-\beta_{n}^{\left(k\right)}\left(x\mid\lambda\right)\\
 & =\sum_{p=1}^{n}\left(\sum_{m=0}^{p-1}\frac{\left(-1\right)^{m+k+1}}{\left(m+1\right)^{k}}\left(m+1\right)!S_{2}\left(k+m+1\right)\right)\binom{n}{p}\left(x\mid\lambda\right)_{n-p}.
\end{align*}

\end{thm}
Now, we note that 
\begin{align}
 & \relphantom{=}{}\frac{\Li_{k}\left(1-e^{-t}\right)}{\left(1+\lambda t\right)^{\frac{1}{\lambda}}-1}\left(1+\lambda t\right)^{\frac{x}{\lambda}}\label{eq:20}\\
 & =\frac{\Li_{k}\left(1-e^{-t}\right)}{\left(1+\lambda t\right)^{\frac{d}{\lambda}}-1}\sum_{a=0}^{d-1}\left(1+\lambda t\right)^{\frac{l+x}{\lambda}}\nonumber \\
 & =\left(\frac{\Li_{k}\left(1-e^{-t}\right)}{t}\right)\frac{1}{d}\sum_{a=0}^{d-1}\frac{dt}{\left(1+\lambda t\right)^{\frac{d}{\lambda}}-1}\left(1+\lambda t\right)^{\frac{l+x}{\lambda}}\nonumber\\
 & =\sum_{l=0}^{\infty}\left(\sum_{p=1}^{l+1}\frac{\left(-1\right)^{p+l+1}}{p^{k}}p!\frac{S_{2}\left(l+1,p\right)}{l+1}\right)\frac{t^{l}}{l!}\nonumber \\
 & \relphantom{=}{}\times\sum_{a=0}^{d-1}\sum_{m=0}^{\infty}\beta_{m}\left(\left.\frac{l+x}{d}\right|\frac{\lambda}{d}\right)d^{m-1}\frac{t^{m}}{m!}\nonumber \\
 & =\sum_{a=0}^{d-1}\left(\sum_{n=0}^{\infty}\left(\sum_{l=0}^{n}\sum_{p=1}^{l+1}\binom{n}{l}\frac{\left(-1\right)^{p+l+1}}{p^{k}}p!\frac{S_{2}\left(l+1,p\right)}{l+1}\beta_{n-l}\left(\left.\frac{l+x}{d}\right|\frac{\lambda}{d}\right)d^{n-l-1}\right)\frac{t^{n}}{n!}\right)\nonumber \\
 & =\sum_{n=0}^{\infty}\left\{ \sum_{a=0}^{d-1}\sum_{l=0}^{n}\sum_{p=1}^{l+1}\binom{n}{l}\frac{\left(-1\right)^{p+l+1}}{p^{k}}p!\frac{S_{2}\left(l+1,p\right)}{l+1}\beta_{n-l}\left(\left.\frac{l+x}{d}\right|\frac{\lambda}{d}\right)d^{n-l-1}\right\} \frac{t^{n}}{n!},\nonumber 
\end{align}
where $d$ is a fixed positive integer. 

On the other hand, 
\begin{align}
 & \relphantom{=}{}\frac{\Li_{k}\left(1-e^{-t}\right)}{\left(1+\lambda t\right)^{\frac{1}{\lambda}}-1}\left(1+\lambda t\right)^{\frac{x}{\lambda}}\label{eq:21}\\
 & =\sum_{n=0}^{\infty}\beta_{n}^{\left(k\right)}\left(x\mid\lambda\right)\frac{t^{n}}{n!}.\nonumber 
\end{align}

Therefore, by (\ref{eq:20}) and (\ref{eq:21}), we obtain the following
theorem.
\begin{thm}
\label{thm:4} For $n\ge0$, $d\in\mathbb{N}$ and $k\in\mathbb{Z}$,
we have 
\begin{align*}
 & \relphantom{=}{}\beta_{n}^{\left(k\right)}\left(x\mid\lambda\right)\\
 & =\sum_{a=0}^{d-1}\sum_{l=0}^{n}\sum_{p=1}^{l+1}\binom{n}{l}\frac{\left(-1\right)^{p+l+1}}{p^{k}}p!\frac{S_{2}\left(l+1,p\right)}{l+1}\beta_{n-l}\left(\left.\frac{l+x}{d}\right|\frac{\lambda}{d}\right)d^{n-l-1}.
\end{align*}

\end{thm}
From (\ref{eq:13}), we can derive the following equation: 
\begin{align}
 & \relphantom{=}{}\sum_{n=0}^{\infty}\beta_{n}^{\left(k\right)}\left(x+y\mid\lambda\right)\frac{t^{n}}{n!}\label{eq:22}\\
 & =\frac{\Li_{k}\left(1-e^{-t}\right)}{\left(1+\lambda t\right)^{\frac{1}{\lambda}}-1}\left(1+\lambda t\right)^{\frac{x+y}{\lambda}}\nonumber \\
 & =\left(\frac{\Li_{k}\left(1-e^{-t}\right)}{\left(1+\lambda t\right)^{\frac{1}{\lambda}}-1}\left(1+t\lambda\right)^{\frac{x}{\lambda}}\right)\left(1+\lambda t\right)^{\frac{y}{\lambda}}\nonumber \\
 & =\left(\sum_{l=0}^{\infty}\beta_{l}^{\left(k\right)}\left(x\mid\lambda\right)\frac{t^{l}}{l!}\right)\left(\sum_{m=0}^{\infty}\left(y\mid\lambda\right)_{m}\frac{t^{m}}{m!}\right)\nonumber \\
 & =\sum_{n=0}^{\infty}\left(\sum_{l=0}^{n}\binom{n}{l}\beta_{l}^{\left(k\right)}\left(x\mid\lambda\right)\left(y\mid\lambda\right)_{n-l}\right)\frac{t^{n}}{n!}.\nonumber 
\end{align}

Therefore, by (\ref{eq:22}), we obtain the following theorem.
\begin{thm}
\label{thm:5} For $n\ge0$, we have 
\[
\beta_{n}^{\left(k\right)}\left(x+y\mid\lambda\right)=\sum_{l=0}^{n}\binom{n}{l}\beta_{l}^{\left(k\right)}\left(x\mid\lambda\right)\left(y\mid\lambda\right)_{n-l}.
\]
\end{thm}
\begin{rem*}
$\,$
\end{rem*}
\begin{align*}
 & \relphantom{=}\frac{d}{dx}\beta_{n}^{\left(k\right)}\left(x\mid\lambda\right)\\
 & =\frac{d}{dx}\sum_{l=0}^{n}\binom{n}{l}\beta_{n-l}^{\left(k\right)}\left(\lambda\right)\left(x\mid\lambda\right)_{l}\\
 & =\sum_{l=0}^{n}\binom{n}{l}\beta_{n-l}^{\left(k\right)}\left(\lambda\right)\sum_{j=0}^{l-1}\frac{1}{x-\lambda j}\prod_{i=0}^{l-1}\left(x-\lambda i\right)\\
 & =\sum_{l=0}^{n}\binom{n}{l}\beta_{n-l}^{\left(k\right)}\left(\lambda\right)\sum_{j=0}^{l-1}\prod_{\substack{i=0\\
i\neq j
}
}^{l-1}\left(x-\lambda i\right).
\end{align*}

{\bf Competing interests}\\
The authors declare that they have no competing interests.\\\\

{\bf Authors' contributions}\\
All authors contributed equally to this work. All authors read and approved the final manuscript.\\\\

\bibliographystyle{amsplain}

\providecommand{\bysame}{\leavevmode\hbox to3em{\hrulefill}\thinspace}
\providecommand{\MR}{\relax\ifhmode\unskip\space\fi MR }
\providecommand{\MRhref}[2]{%
  \href{http://www.ams.org/mathscinet-getitem?mr=#1}{#2}
}
\providecommand{\href}[2]{#2}

\end{document}